\documentclass[12pt]{amsart}

\usepackage{amssymb, amsfonts, amsthm}
\usepackage{graphicx}
\usepackage{color}

\numberwithin{equation}{section}

\newtheorem{thm}{Theorem}[section]
\newtheorem{cor}[thm]{Corollary}
\newtheorem{lem}[thm]{Lemma}

\theoremstyle{definition}
\newtheorem{defn}[thm]{Definition}
\newtheorem{rem}[thm]{Remark}

 \numberwithin{equation}{section}
\newtheorem{conj}[thm]{Conjecture}
\def\be{\begin{equation}}
\def\ee{\end{equation}}
\newcommand{\bthm}{\begin{thm}}
\newcommand{\ethm}{\end{thm}}
\newcommand{\bcor}{\begin{cor}}
\newcommand{\ecor}{\end{cor}}
\newcommand{\bcon}{\begin{conj}}
\newcommand{\econ}{\end{conj}}
\newcommand{\blem}{\begin{lem}}
\newcommand{\elem}{\end{lem}}
\newcommand{\bdefn}{\begin{defn}}
\newcommand{\edefn}{\end{defn}}

\newcommand{\beq}{\begin{eqnarray}}
\newcommand{\beqq}{\begin{eqnarray*}}
\newcommand{\eeq}{\end{eqnarray}}
\newcommand{\eeqq}{\end{eqnarray*}}
\newcommand{\bpf}{\begin{proof}}
\newcommand{\epf}{\end{proof}}

\newcommand{\brem}{\begin{rem}}
\newcommand{\erem}{\end{rem}}

\newcommand{\bconj}{\begin{conj}}
\newcommand{\econj}{\end{conj}}

\newcommand{\A}{{\mathcal A}}

\newcommand{\T}{{\mathcal T}}

\newcommand{\es}{{\mathcal S}}

\newcommand{\IC}{{\mathbb C}}

\newcommand{\D}{{\mathbb D}}

\newcommand{\N}{{\mathbb N}}

\newcounter{minutes}
\divide\time by 60
\newcounter{hours}
\multiply\time by 60 \addtocounter{minutes}{-\time}
\begin{document}

\title[Logarithmic coefficients for Janowski type $(j,k)$-symmetric starlike functions]
{A correct proof of logarithmic coefficients for Janowski type $(j,k)$-symmetric starlike functions
}

\author[Navneet Lal Sharma]{Navneet Lal Sharma$^{1,2}$}
\address{Navneet Lal Sharma,
$^1 $Department of Mathematics, 
National Rail and Transportation Institute (NRTI) Vadodara.
NAIR Campus, Lalbaug, Vadodara-390004, Gujarat. India }
\email{sharma.navneet23@gmail.com \&  nlsharma@nrti.edu.in}

\address{$^2$ School of Mathematical Science, Universiti Sains of Malaysia, Penang, Malaysia. 
}


\begin{abstract}
Recently, authors~\cite{HG22} studied the logarithmic coefficient bounds for class of the Janowski type $(j,k)$-symmetric starlike functions $\mathcal{ST}_{[j,k]}(A,B)$ in ({\em Rev. Real Acad. Cienc. Exactas Fis. Nat. Ser. A-Mat.} (2022)). We pointed out that the proof of Theorem~3 in \cite{HG22} is incorrect. 
In this article, we present the correct proof of Theorem~3. 
In addition, we also obtain some new results related to the logarithmic coefficient inequalities for the class $\mathcal{ST}_{[j,k]}(A,B)$.

\end{abstract}

\maketitle \pagestyle{myheadings} \markboth{ N. L. Sharma}{Logarithmic coefficient for Janowski type $(j,k)$-symmetric starlike functions}


\begin{center}
\texttt{File:~\jobname .tex, printed: \number\day-\number\month-\number\year,
\thehours.\ifnum\theminutes<10{0}\fi\theminutes}
\end{center}

\thispagestyle{empty}

\section{Preliminaries} 
Let $\mathcal{H}$ be the class of all analytic functions in the unit disk $\D:=\{z\in\mathbb{C}:\, |z|<1 \}$
and  $\A$ the subclass of $\mathcal{H}$ consisting of functions $f$ with normalization $f(0)=0=f'(0)=1$. Any function $f\in\A$ has the Taylor series expansion 
$ f(z)=z+\sum_{n=2}^\infty a_n z^n,\quad z\in\D. $
Let $\es$ denote the subclass of $\A$ consisting of functions that are univalent in $\D$.
The class $\es$ is invariant under the $k$th root transformation. 
For each function $f\in \A$, the function $(f(z^k))^{1/k}$ for $k\in \mathbb{N}$ is univalent and maps the unit disc $\D$ into a region with $k$-fold symmetry. Every function with $k$-fold symmetric has the following series representation:
$$ f_k(z)=z+\sum_{n=1}^\infty a_{nk+1} z^{nk+1},\quad z\in\D,\,k\in \mathbb{N}. 
$$
Denote by $\Omega$ the class of $k$-fold symmetric. A function $f:\Omega \rightarrow \mathbb{C}$ is called $(j,k)$ symmetrical if
$$ f(\epsilon z) =\epsilon ^j f(z) \quad \left(\mbox{for all } z\in\Omega; \quad \epsilon =e^{2\pi i /k}; \quad j=0,1,2,\cdots,k-1 \mbox{ for } k\ge 2 \right).
$$
In this paper, we consider the domain $\Omega$ as $\D$ and denote by $\es ^{(j,k)}$ the class of all ${(j,k)}$-symmetric functions in $\D$. 
Initially, Liczberski and Polubinski introduced and study the ${(j,k)}$-symmetrical function class $\es ^{(j,k)}$. 
If function $f\in \mathcal{A}$ is $(j,k)$-symmetric, then $f_{j,k}$ has the following series expansion:
\beq\label{eq1}
f_{j,k}(z)=z+\sum_{n=1}^\infty a_{nk+j} z^{nk+j} \quad \mbox{ for } z\in\D. 
\eeq
For $-1\leq B\leq 0$ and $A\in  \IC$, $A\neq B$, let $\es \T_{[j,k]}(A,B)$ denote the class of ${(j,k)}$-symmetric functions of the form (\ref{eq1})
which satisfy the subordination relation
$$\es \T_{[j,k]}(A,B):=\left \{f\in \A: \frac{zf'_{j,k}(z)(z)}{f_{j,k}(z)}\prec \frac{1+Az^{j+k-1}}{1+Bz^{j+k-1}},
\quad   z\in \D   \right\}
$$
where $j=0,1,2,\cdots,k-1,\, k\in \N$ and the symbol $\prec$ denotes the usual subordination.
In particular, $\es^*(A,B):=\es \T_{[1,1]}(A,B)$ is the class of Janowski type starlike functions and 
$\es^*_{k}(A,B):=\es \T_{[1,k]}(A,B)$ is the class of Janowski type $k$-fold symmetric starlike functions.
For other suitable choice of the parameters $j,k,A$ and $B$, we can obtain different subclasses studied by various authors, see~\cite{HG22}.
The function $K^{(j,k)}_{A,B}$ defined by
\begin{equation}\label{eq2}
\displaystyle K^{(j,k)}_{A,B}(z):=
\left \{
\begin{array}{ll}
ze^{\frac{A\,z^{j+k-1}}{j+k-1}} & \mbox{ for } B=0\\
z(1+Bz^{j+k-1})^{\frac{A-B}{(j+k-1)B}} & \mbox{ for } B\neq 0
     \end{array}, \right.
\end{equation}
is in $\es \T_{[j,k]}(A,B)$ and plays the role of extremal function for the class $\es \T_{[j,k]}(A,B)$. 

In addition to analytic and geometric considerations, logarithmic restrictions and special exponentiation techniques are frequently helpful in the theory of univalent functions. Milin examined the effect of transferring the properties of logarithmic coefficients to those of Taylor coefficients of univalent functions themselves or to their powers during the 1960s, and their role in the theory of univalent functions. Milin's inequalities received a lot of interest because their validity would imply the truth of the Robertson and Bieberbach conjectures. 
Louis de Branges proved these inequalities and also settled the Bieberbach conjecture for the class $\es.$ The proof which settles the
Bieberbach conjecture relied not on the coefficients $\{a_n\}$ of $f$ but rather the logarithmic coefficients of $f\in\es$. For more details on the logarithmic coefficients, see~\cite{PSW1,PSW2}.

In the case of  ${(j,k)}$-symmetric functions, the logarithmic coefficients $d_{n(j+k-1)}$ of  ${(j,k)}$-symmetric functions defined by the formula 
\begin{equation}\label{eq3}
\log\left( \frac{f_{j,k}(z)}{z}\right) =2\sum_{n=1}^{\infty} d_{n(j+k-1)} z^{n(j+k-1)}, \quad z\in \D.
\end{equation}
In~\cite{HG22}, Srivastava {\it et al.} recently established the estimate for the $n$th logarithmic coefficients of functions in the class $\es \T_{[j,k]}(A,B)$. \\

\noindent
{\bf Theorem~A ~\cite[Theorem~3]{HG22}.} 
Suppose $f_{j,k} \in \es \T_{[j,k]}(A,B)$ for $A\in\IC $,\,  $-1\le B\le 0$ and $A\neq B$. Then the corresponding logarithmic coefficient satisfies the following inequality:
\begin{equation}\label{eq4}
\sum_{n=1}^{\infty}| d_{n(j+k-1)}|^{2} \leq \left(\frac{|A-B|}{2(j+k-1)B}\right)^2 {\rm Li\,}_{2}(B^2),
\end{equation}
where 
$${\rm Li\,}_{v}(z)=\sum_{n=1}^{\infty}\frac{z^n}{n^v} \quad \quad (v\in\IC \mbox{ and } z\in\D;\, {\rm Re}(v)>1 \mbox{ and } |z|=1)
$$ denotes the Polylogarithm function (or de Jonqui\'{e}re's function) of order $v$. The estimate in (\ref{eq4}) is the best possible..

We pointed out that the proof of theorem~A is incorrect. 
In this paper, we present the correct proof of Theorem~A. 
In addition, we also determine some results related the logarithmic coefficient inequalities for the class $\mathcal{ST}_{[j,k]}(A,B)$.

\section{Logarithmic Coefficients for the Class $\mathcal{ST}_{[j,k]}(A,B)$} 

In this section, first we provide a correct proof of Theorem~A.

\subsection*{Proof of Theorem~A}
Suppose $f_{j,k} \in \es \T_{[j,k]}(A,B)$. Then by the definition of $\es \T_{[j,k]}(A,B)$, we get
\begin{equation}\label{eqA.1}    
\frac{zf'_{j,k}(z)}{f_{j,k}(z)}\prec \frac{1+Az^{j+k-1}}{1+Bz^{j+k-1}}, \quad   z\in \D . 
\end{equation}
Let $g(z):=z/f_{j,k}$ which is a non-vanishing analytic function in $\D$ with $g(0)=1$ and satisfies the relation
\begin{equation}\label{eqA.2}   
\frac{zg'(z)}{ g(z)} = 1- \frac{zf'_{j,k}(z)}{f_{j,k}(z)}. 
\end{equation}
From the relations (\ref{eqA.1}) and (\ref{eqA.2}), we get
$$ \frac{zg'(z)}{g(z)}\prec \frac{-(A-B)z^{j+k-1}}{1+B z^{j+k-1}}=:G(z), \quad z\in \D.
$$
Since $G$ is convex in $\D$ and $G(0)=0$, then by the well-known subordination result~\cite[Corollary~3.1d.1, p.~76]{MM2000}, we obtain
\begin{equation}\label{eqA.3}
g(z)=\frac{z}{f_{j,k}(z)} \prec Q^{(j,k)}_{A,B}(z)=\exp \left(\int_0^{z}\frac{G(t)}{t} dt \right).
\end{equation}
It is a simple calculation to compute that
$$ \displaystyle Q^{(j,k)}_{A,B}(z)=
 \left \{
 \begin{array}{ll}
 e^{-\frac{Az^{j+k-1}}{j+k-1}} & \mbox{ for } B=0\\
 (1+Bz^{j+k-1})^{\frac{-(A-B)}{({j+k-1})B}} & \mbox{ for } B\neq 0.
 \end{array} \right.
$$
We can rewrite the relation (\ref{eqA.3}) as
$$ \frac{f_{j,k}(z)}{z}\prec \frac{1}{Q^{(j,k)}_{A,B}(z)},
$$
which, in terms of the logarithmic coefficients $ d_{n(j+k-1)}$ of $f_{j,k}$ defined by (\ref{eq3}), is equivalent to
$$ 2\sum_{n=1}^{\infty} d_{n(j+k-1)}z^{n(j+k-1)} \prec
\left \{
\begin{array}{ll}
\displaystyle
\left(\frac{A-B}{(j+k-1)B}\right)\sum_{n=1}^{\infty} (-1)^{n-1}\frac{B^n}{n}z^{n(j+k-1)}
& \mbox{ for } B\neq 0\\[6mm]
\frac{Az^{j+k-1}}{j+k-1}
& \mbox{ for } B= 0.
\end{array}\right.
$$
Then by using Rogosinski's Theorem~\cite[Theorem~6.3, p.~192]{Dur83}\, (see~\cite{Rog43}), we obtain that
$$  4\sum_{n=1}^{k} |d_{n(j+k-1)}|^2  \le \frac{1}{(j+k-1)^2}\left|\frac{A-B}{B}\right|^2\sum_{n=1}^{k} \frac{|B|^{2n}}{n^2}
 \quad \mbox{ for }B\neq 0.
 $$   

Applying $k \rightarrow \infty$, we get
\begin{align*} \sum_{n=1}^{\infty} |d_{n(j+k-1)}|^2 & \le \left(\frac{|A-B|}{2(j+k-1)B}\right)^2\sum_{n=1}^{\infty} \frac{B^{2n}}{n^2}\\
& = \frac{|A-B|^2}{4(j+k-1)^2}\frac{{\rm Li\,}_{2}(B^2)}{B^2},
\end{align*}
where ${\rm Li\,}_{2}(x)=\sum_{n=1}^{\infty}\frac{x^n}{n^2}$. For $x=0$, we let $\frac{{\rm Li\,}_{2}(x)}{x}$ as the limit value 1.
This proves the desired assertion (\ref{eq4}).

The equality holds for the function $K^{(j,k)}_{A,B}$  defined by (\ref{eq2}).
Indeed, for the function $K^{(j,k)}_{A,B}$, we have
\begin{align*}
\displaystyle \log \left( \frac{K^{(j,k)}_{A,B}(z)}{z} \right) & =
\left \{
\begin{array}{ll}
 \displaystyle \left (\frac{A-B}{(j+k-1)B}\right )\log (1+Bz^{j+k-1}) & \mbox{ for } B\neq 0\\[6mm]
\frac{Az^{j+k-1}}{j+k-1} & \mbox{ for } B=0
 \end{array} \right. \\[4mm]
 & =: 2\sum_{n=1}^{\infty} d_{n(j+k-1)}(K^{(j,k)}_{A,B})\,z^{n(j+k-1)},
\end{align*}
where
\begin{align}\label{eqA.4}
\displaystyle 2\,d_{n(j+k-1)}(K^{(j,k)}_{A,B}) =
\left \{
\begin{array}{ll}
\displaystyle 
(-1)^{n-1} \left (\frac{A-B}{(j+k-1)B}\right )\frac{B^n}{n}& \mbox{ for } B\neq 0 \\[2mm]
\displaystyle  \frac{A}{j+k-1} & \mbox{ for } B=0. 
 \end{array} \right. \\ \nonumber
\end{align}
This completes the proof of Theorem~A.
\hfill$\Box$
\bcor\label{cor1}
 If $f_{j,k} \in \es \T_{[j,k]}(A,0)$ for $A\in \IC\backslash \{0\}$, then the logarithmic coefficients of $f_{j,k}$ satisfy 
 $$ \sum_{n=1}^{\infty}|d_{n(j+k-1)}|^2 \le \frac{|A|^2}{4(j+k-1)^2}.
 $$
The inequality is sharp for the functions $K^{(j,k)}_{A,0}=ze^{\frac{A\,z^{j+k-1}}{j+k-1}}.$\\
In particular, for $f_{k} \in \es^*_{k}(A,0)$, one has the sharp inequality
$ \sum_{n=1}^{\infty}|d_{nk}|^2 \le |A|^2/4k^2 .$ 
\ecor
\bcor\label{cor2}
 Let $f_{j,k} \in \es \T_{[j,k]}(A,-A)$ for $0<A \le 1$. Then we get 
 $$ \sum_{n=1}^{\infty}|d_{n(j+k-1)}|^2 \le \frac{1}{(j+k-1)^2} {\rm Li_2}(A^2).
 $$
The inequality is sharp for the functions $K^{(j,k)}_{A}=z(1-Az^{j+k-1})^{\frac{-2}{(j+k-1)}}.$\\
In particular, if $f_{k} \in \es^*_{k}(A,-A)$, then we obtain the sharp inequality
$ \sum_{n=1}^{\infty}|d_{nk}|^2 \le {\rm Li_2}(A^2)/k^2.$
\ecor

Next we state the logarithmic coefficients inequality for the class $\es \T_{[j,k]}(A,B)$.
\bthm\label{thm2}
For $A\in\IC $,\,  $-1\le B\le 0$ and $A\neq B$, the logarithmic coefficients of $f_{j,k} \in \es \T_{[j,k]}(A,B)$ satisfy the inequality
$$ \sum_{n=1}^{\infty} n^2|d_{n(j+k-1)}|^{2} \leq \frac{1}{4(j+k-1)^2} \frac{|A-B|^2}{1-B^2} \quad \mbox{ for } B\neq -1.
$$
The inequality is sharp for the function $K^{(j,k)}_{A,B}$ defined by (\ref{eq2}).
\ethm

\bpf
Suppose $f_{j,k} \in \es \T_{[j,k]}(A,B)$, then from the relation (\ref{eqA.1}), we obtain
$$ z \frac{d}{dz}\left[ \log \left( \frac{f_{j,k}(z)}{z} \right)\right]=
\frac{zf'_{j,k}(z)}{f_{j,k}(z)} - 1 \prec \frac{(A-B)z^{j+k-1}}{1+B z^{j+k-1}}, \quad   z\in \D
$$
which, in terms of the logarithmic coefficients $ d_{n(j+k-1)}$ of $f_{j,k}$ defined by (\ref{eq3}), is equivalent to
$$ 2(j+k-1)\sum_{n=1}^{\infty}n d_{n(j+k-1)}z^{n(j+k-1)} \prec
\left \{
\begin{array}{ll}
\displaystyle
\left(\frac{A-B}{B}\right)\sum_{n=1}^{\infty} (-1)^{n-1} B^n z^{n(j+k-1)} & \mbox{ for } B\neq 0\\[6mm]
Az^{j+k-1} & \mbox{ for } B= 0.
\end{array}\right.
$$
By using Rogosinski’s result~\cite[Theorem~II (i)]{Rog43}, we get the inequality
\beq\label{eq5} 
4(j+k-1)^2 \sum_{n=1}^{k} n^2 \, |d_{n(j+k-1)}|^2  \le \left|\frac{A-B}{B}\right|^2\sum_{n=1}^{k} |B|^{2n}  \quad \mbox{ for }B\neq -1.
 \eeq 
Applying $k \rightarrow \infty$, we get
\begin{align*} 
\sum_{n=1}^{\infty} n^2 \, |d_{n(j+k-1)}|^2 & \le \frac{1}{4(j+k-1)^2}\left(\frac{|A-B|}{B}\right)^2\sum_{n=1}^{\infty}B^{2n} \quad \mbox{ for } B\neq -1\\
& = \frac{1}{4(j+k-1)^2} \frac{|A-B|^2}{1-B^2}  \quad \mbox{ for }B\neq -1,
\end{align*}
which leads us to the required result.
The equality holds for the function $K^{(j,k)}_{A,B}$ and logarithmic coefficients $ d_{n(j+k-1)}$ of $K^{(j,k)}_{A,B}$ is defined by (\ref{eqA.4}).
\epf

If we choose $j=1$ in Theorem~\ref{thm2}, then we obtain the following result for $\es^*_k(A,B)$.
\bcor
If  $f_{k} \in \es^*_k(A,B)$ for $A\in\IC $,\,  $-1\le B\le 0$ and $A\neq B$, then the logarithmic coefficients of $f_{k} $ satisfy the inequality
$$ \sum_{n=1}^{\infty} n^2|d_{nk}|^{2} \leq \frac{1}{4k^2} \frac{|A-B|^2}{1-B^2} \quad \mbox{ for } B\neq -1.
$$
The inequality is sharp for the function $K^{(k)}_{A,B}$ defined by (\ref{eq2}).
\ecor

Our next result, establishes an inequality of the type~\cite[Theorem~1.1]{Rot07} for the class $\es \T_{[j,k]}(A,B)$.
\bthm\label{thm3}
Let $f_{j,k} \in \es \T_{[j,k]}(A,B)$ for $A\in\IC $,\,  $-1\le B\le 0$ and $A\neq B$, and let $t\le 2$. Then the logarithmic coefficients of $f_{j,k}$ satisfy the inequality
$$ \sum_{n=1}^{\infty} (n+1)^t \, |d_{n(j+k-1)}|^{2} \leq \left( \frac{|A-B|}{2(j+k-1)B} \right)^2 \sum_{n=1}^{\infty} \frac{(n+1)^t}{n^2} B^{2n}.
$$
\ethm

\bpf
Suppose  $f_{j,k} \in \es \T_{[j,k]}(A,B)$. We recall from relation (\ref{eq5}), that for $k\in\N$ the inequalities
\beq\label{eq6} 
\sum_{n=1}^{k} n^2 |d_{n(j+k-1)}|^2  \le  \left(\frac{|A-B|}{2(j+k-1)B}\right)^2 \sum_{n=1}^{k} B^{2n} : = H(A,B)  \sum_{n=1}^{k} B^{2n}
\eeq
are valid. 
We now take inequality (\ref{eq6}) into account for  $k=1,2,3,\cdots,N,$ and multiply the $k$th inequality by the factor
$$ \frac{(k+1)^t}{k^2} - \frac{(k+2)^t}{(k+1)^2} >0 \quad \mbox{ for } k= 1,2,3,\cdots ,N-1,
$$
and multiply the $N$-th inequality by the factor $(N+1)^t / N^2$ for $k=N$.
Now, we obtain in the left-hand side of the inequality of (\ref{eq6}) after summing all of these modified inequalities
\begin{align*}
 \sum_{k=1}^{N-1} & \left[ \left(\frac{(k+1)^t}{k^2} - \frac{(k+2)^t}{(k+1)^2} \right) \sum_{n=1}^k n^2 \, |d_{n(j+k-1)}|^2  \right]
 + \frac{(N+1)^t}{N^2} \sum_{n=1}^N n^2 \, |d_{n(j+k-1)}|^2  \\
 & =\sum_{n=1}^{N-1} n^2 \, |d_{n(j+k-1)}|^2  \frac{(n+1)^t}{n^2} + (N+1)^t\,|d_{N(j+k-1)}|^2  \\
 & = \sum_{n=1}^{N} (n+1)^t \, |d_{n(j+k-1)}|^2 
\end{align*}
furthermore, on the right side of the inequality of (\ref{eq6}), we find
\begin{align*}
 & H(A,B)\left(\sum_{k=1}^{N-1}  \left[ \left(\frac{(k+1)^t}{k^2} - \frac{(k+2)^t}{(k+1)^2} \right) \sum_{n=1}^k  B^{2n}  \right]
 + \frac{(N+1)^t}{N^2} \sum_{n=1}^N  B^{2n} \right) \\
 & \hspace{3cm} =H(A,B) \left(\sum_{n=1}^{N-1} \frac{(n+1)^t}{n^2} B^{2n} + \frac{(N+1)^t}{N^2} B^{2N} \right)  \\
 & \hspace{3cm} = H(A,B)\sum_{n=1}^{N} \frac{(n+1)^t}{n^2} B^{2n} \quad \quad (H(A,B) \mbox{ defined in } (\ref{eq6})). 
\end{align*}
As a result, we get the following inequality
$$ \sum_{n=1}^{N} (n+1)^t \, |d_{n(j+k-1)}|^2 \le  \left(\frac{|A-B|}{2(j+k-1)B}\right)^2  \sum_{n=1}^{N} \frac{(n+1)^t}{n^2} B^{2n}. 
$$
Finally, letting $n\rightarrow \infty$, we obtain
$$ \sum_{n=1}^{\infty} (n+1)^t \, |d_{n(j+k-1)}|^2 \le  \left(\frac{|A-B|}{2(j+k-1)B}\right)^2  \sum_{n=1}^{\infty} \frac{(n+1)^t}{n^2} B^{2n}. 
$$
\epf

\brem
Here is an alternate approach to prove the inequality (\ref{eq4}). If we take $t=0$ in Theorem~\ref{thm3}, then we obtain
 $$ \sum_{n=1}^{\infty} |d_{n(j+k-1)}|^2 \le  \left(\frac{|A-B|}{2(j+k-1)B}\right)^2  \sum_{n=1}^{\infty}\frac{1}{n^2} B^{2n}
= \left(\frac{|A-B|}{2(j+k-1)}\right)^2  \frac{{\rm Li\,}_{2}(B^2)}{B^2}. 
$$
\erem
\brem
If we choose $j=1=k$ in Theorems~\ref{thm2},~\ref{thm3} and in Corollaries~\ref{cor1},~\ref{cor2}, then we obtain the results of 
Ponnusamy {\em et al.}~\cite{PSW1}.
\erem
For the choice $j=1$, Theorem~\ref{thm3} reduces to the following logarithmic coefficients inequality for the class $\es^*_k(A,B)$.
\bcor
Let $f_{k} \in \es^*_k(A,B)$ for $A\in\IC $,\,  $-1\le B\le 0$ and $A\neq B$, and let $t\le 2$. Then we obtain
$$ \sum_{n=1}^{\infty} (n+1)^t \, |d_{nk}|^{2} \leq \frac{|A-B|^2}{4k^2B^2} \sum_{n=1}^{\infty} \frac{(n+1)^2}{n^2} B^{2n}.
$$
\ecor


\bigskip
\noindent
{\bf Acknowledgements.} 
Science and Engineering Research Board, Department of Science and Technology, India is supporting the author's work through the "SERB International Research Experience (SIR/2022/000764)" scheme. This work was done while the author was at Universiti Sains Malaysia, under SIRE scheme.

\end{document}